%Paper: alg-geom/9408005
%From: "Dr A.D. King" <aking@liverpool.ac.uk>
%Date: Wed, 10 Aug 1994 17:45:25 +0100 (BST)

%%%%%%%%%%%%%%%%%%%%%%%%%%%%%%%%%%%%%%%%%%%%%%%%%%%%%%%%%%
%%%
%%%   Moduli of Brill-Noether pairs on algebraic curves
%%%
%%%   A.D. King, P.E. Newstead
%%%
%%%   Plain TeX
%%%
%%%   Last revised: 10 Aug 94
%%%
%%%%%%%%%%%%%%%%%%%%%%%%%%%%%%%%%%%%%%%%%%%%%%%%%%%%%%%%%%

%%%
%%% style
%%%

\magnification=1200
\parindent=1em
\font \smcap=cmcsc10
\font \titlefont=cmbx12
\font \authorfont=cmcsc10
\def \and{{\rm and}\ }

\outer\def \beginsection#1\par{\vskip0pt plus.3\vsize\penalty-200
  \vskip0pt plus-.3\vsize\bigskip\vskip\parskip
  \message{#1}\leftline{\bf#1}\nobreak\smallskip}
\outer\def \proclaim #1. #2\par{\medbreak
  {\bf#1.\enspace}{\sl#2}\par\medbreak}
\outer\def \proof{\medbreak{\bf Proof.}\enspace}
\outer\def \endproof{\medbreak}
\outer\def \say #1.
  {\medbreak{\bf #1.}\enspace}
\outer\def \endsay{\medbreak}
\outer\def \para(#1){\medbreak\noindent{\bf#1.}\enspace}
\def \cite#1{{\bf [#1]}}

%%%
%%% local macros
%%%

\def \cC{{\cal C}}
\def \cG{{\cal G}}
\def \cL{{\cal L}}
\def \cQ{{\cal Q}}
\def \cO{{\cal O}}
\def \cR{{\cal R}}
\def \PP{{\bf P}}
\def \al{\alpha}
\def \eps{\varepsilon}
\def \th{\theta}

\def \om{\omega}
\def \kap{\kappa}
\def \lam{\lambda}
\def \Lam{\Lambda}

\def \RR{{\bf R}}
\def \ZZ{{\bf Z}}

\def \LmE{\Lambda}
\def \LmF{\Pi}
\def \Lm{L}
\def \Rss{{\cR^{ss}}}

\def \sub{\subset}
\def \tsr{\otimes}
\def \into{\hookrightarrow}
\def \isom{\cong}

\def \LRA#1{\>{\buildrel {#1}\over\longrightarrow}\>}
\def \rest#1{|_{#1}}

\def \rk{\mathop{\rm rk}\nolimits}
\def \Srk{\mathop{\hbox{\rm $a$-rg}}}
\def \dg{\mathop{\rm deg}\nolimits}
\def \mult{\mathop{\rm m}\nolimits}
\def \mX{{\rm m}_X}

\def \Hom{\mathop{\rm Hom}\nolimits}
\def \im{\mathop{\rm im}\nolimits}
\def \max{\mathop{\rm max}}
\def \quot{\mathop{\hbox{\smcap Quot}}\nolimits}
\def \grass{\mathop{\hbox{\smcap Gr}}\nolimits}
\def \gr{\mathop{\rm gr}\nolimits}

%%%
%%% title
%%%

\centerline{\titlefont Moduli of Brill-Noether pairs on algebraic curves}
\medskip
\centerline{\authorfont A.D. King \and  P.E. Newstead
\footnote*{{\sevenrm This work was carried out in the framework of the
VBAC Research Group of Europroj.
The first author is supported by a grant from the S.E.R.C.}}}
\bigskip
\centerline{10 August 1994}
\bigskip
\centerline{Dept. of Pure Maths, The University of Liverpool,
Liverpool L69 3BX, U.K.}
\centerline{e-mail:\  aking@liverpool.ac.uk,\  newstead@liverpool.ac.uk}
\bigskip

%%%
%%% main document
%%%

\beginsection 1. Introduction

In this paper we shall construct moduli spaces of pairs consisting of
a torsion-free sheaf $E$ over an algebraic curve $X$
(possibly singular, but always of pure dimension)
and a vector subspace $\LmE$ of its space of sections $H^0(E)$.
Such pairs are the common generalisation of the classical notion of
a linear system, for which $E$ is (the sheaf of sections of)
a line bundle, and of Bradlow pairs \cite{3}, for which one may
take $\LmE$ to be 1-dimensional --- although strictly speaking
$\LmE$ should be replaced by a single non-zero section $\phi$.
These more general pairs have also been studied by Bertram \cite{2},
Raghavendra and Vishwanath \cite{13} (under the name of ``pairs'')
and by Le Potier \cite{8} (under the name of ``coherent systems''),
who works over varieties (or schemes) of arbitrary dimension.
In this paper, we shall use the term ``Brill-Noether pair'',
to emphasise the role that we hope they will play in higher rank
Brill-Noether theory.

There are three numerical invariants which make up the `type'
of a Brill-Noether pair:
the rank $r$ and degree $d$ of the sheaf $E$,
which can be defined even on a singular curve provided it is polarised
(see Definition 2.1), and the dimension $l$ of $\LmE$.
In the case $r=1$ (and $X$ smooth) there is a well-known parameter
space for all linear systems --- usually denoted $G^r_d$
(e.g. \cite{1} Chap IV),
where $r+1$ is the dimension of the space of sections!
In higher rank, one must introduce a notion of semistability,
which depends on a single parameter $\al$ (see Definition 2.3.2).
The same definition is used in \cite{13}, while in \cite{3}
an equivalent definition is used involving $\tau={d+\al l\over r}$.
The result we shall prove is the following.

\proclaim Theorem 1.
Let $X$ be a polarised algebraic curve and
$\al$ any positive rational number.
For each $(r,d,l)$ there exists a projective scheme
$G_\al(r,d,l)$ which is a coarse moduli space for families of
$\al$-semistable Brill-Noether pairs of type $(r,d,l)$.
The (closed) points of $G_\al(r,d,l)$ are in one-one correspondence
with S-equivalence classes of pairs.

For the definition of a family of Brill-Noether pairs see Section 2.5.
As usual two pairs are S-equivalent if they have the same composition factors
in the appropriate category of semistable pairs (c.f. Corollary 2.5.1).
The proof of Theorem 1 is given in Section 3.

When $X$ is smooth, the moduli spaces $G_\al(r,d,l)$ have
already been constructed in different ways and in various cases.
A gauge theoretic construction gives $G$ as a K\"ahler manifold
instead of a projective variety;
some cases are covered by \cite{3}, while the general case is described
in \cite{4}.
An algebraic construction for moduli of Bradlow pairs was
given by Huybrechts and Lehn \cite{6}, regarding the dual
of the section $E^*\to\cO_X$ as the `pair'.
An algebraic construction for general $(r,d,l)$, but with
$\al$ sufficiently small is given in \cite{13}
(and also in \cite{2}, but just for large $d$.)
The case when $\al$ is sufficiently large is covered by \cite{8},
which also allows $X$ to be singular.
More recently, in \cite{9}, Le Potier has observed that one may
also introduce parameters into his construction.

The spaces $G=G_\al$ with $\al$ sufficiently small are the most important
from the point of view of Brill-Noether theory, because they are the spaces
that map to the Brill-Noether loci $W$ in the moduli spaces of
semistable bundles (or, on singular curves, torsion-free sheaves),
i.e. the loci consisting of bundles whose space of sections has
at least a given dimension.
However, inspired by the success of Thaddeus' work \cite{16},
one can expect to obtain useful information about the space
$G$ by studying how $G_\al$ varies as the parameter $\al$ varies.
It is therefore important to be able to construct all such
moduli spaces.
In addition, degeneration arguments (and a desire for completeness)
require one to be able to construct such spaces over singular
curves.
Using the approach developed by Simpson \cite{15},
this can be done without much additional complication.

%%%
%%%
%%%

\beginsection 2. Definitions and Preliminary Results

Throughout this paper $X$ is a polarised algebraic curve
defined over an algebraically closed field $k$.
Here ``curve'' means a connected scheme of pure dimension 1,
while ``polarised'' means that $X$ is equipped with
an ample invertible sheaf $\cO_X(1)$.
(Note also that ``sheaf'' always means a coherent sheaf
of $\cO_X$-modules.)
Without loss of generality, we shall assume that $\cO_X(1)$
is `primitive', i.e. not a power ($>1$) of another invertible
sheaf.
Let $N_X$ be the smallest integer such that $\cO_X(N_X)$ is
very ample.
In addition, let $-\beta_X$ be the minimum slope amongst
all quotient sheaves of $\cO_X$.
(On a non-integral curve $\cO_X$ may not be stable or even
semistable.)
It is easy to see then that $-\beta_X$ is the minimum slope of
any sheaf that is generated by its sections.

\para(2.1)
For any sheaf $E$, the `multiplicity' $\mult(E)$
(c.f. \cite{8} Def. 2.1) is defined as the leading coefficient
of the Hilbert polynomial
$$
\chi\bigl(E(n)\bigr)=\mult(E)n+\chi(E).
$$
We then define the `multiplicity of $X$' $\mX=\mult(\cO_X)$.
Note that this depends on the polarisation and the fact that it is
primitive.
In particular, $\mX=1$ if and only if $X$ is integral
and then the primitive $\cO_X(1)$ has degree 1 in the usual sense.

One may work with $\mult(E)$ and $\chi(E)$ as the two numerical
invariants of a sheaf $E$,
but we prefer to work with two others which directly generalise
the rank and degree in the integral (in particular, smooth) case.
As a bonus, these invariants do not depend on the
polarisation being primitive, however one must be careful because
they are not necessarily integers, but only integer multiples of
$1/\mX$.

\say Definition 2.1.
For a sheaf $E$, the \sl rank \rm $\rk E$ and \sl degree \rm $\dg E$
are defined from the Hilbert polynomial as follows
$$
\chi\bigl(E(n)\bigr)=(\rk E) \chi\bigl(\cO_X(n)\bigr) + \dg E.
$$
As usual, the \sl slope \rm $\mu(E)=\deg E/\rk E$.
\endsay

When $X$ is an integral curve, then one readily checks that
$\rk E$ and $\dg E$ are the usual rank and degree.
When $X$ is a reduced curve, then $\rk E$ is equal to Seshadri's
$\Srk(E)$ (\cite{14} p.153),
where $a$ is determined by the degrees of $\cO_X(1)$ along the components.
In general, $\rk E=\mult(E)/\mX$.

Note that, if $E$ is locally free and actually locally isomorphic
to $\cO_X^r$, then $\rk E=r$, as one would expect.
Since this is an integer, $\deg E$ is also an integer in this case.

\para(2.2)
One may also observe that $\rk E=0$ if and only if $E$ is a `torsion sheaf',
i.e. a sheaf with 0-dimensional support.
We will use the term ``torsion-free'' for a sheaf with
no non-trivial torsion subsheaves.
Other terms used for such sheaves are ``of pure dimension 1''
(\cite{15}, \cite{8})
and ``of depth 1'' (\cite{14}).
Because $X$ has pure dimension, the structure sheaf $\cO_X$
and the dualising sheaf $\om_X$ are both torsion-free.

We can now give a crude bound on the dimension of the spaces of sections
of any torsion-free sheaf.

\proclaim Lemma 2.2.
Let $E$ be a torsion-free sheaf such that every non-zero subsheaf $F$ of $E$
has slope $\mu(F)\leq b$.
If $b +\beta_X <0$, then $H^0(E) =0$.
Otherwise, if $b +\beta_X \geq0$, then
$$
\dim H^0(E) \leq \rk E \left( b +\beta_X +\mX N_X \right)
$$

\proof
If every subsheaf of $E$ has slope less than $-\beta_X$, then
the subsheaf generated by the sections of $E$ must be zero, i.e. $H^0(E)=0$.
Otherwise, let $n$ be the positive integer $n$ for which
$n-1 \leq{b +\beta_X \over\mX N_X} <n$.
Then every non-zero subsheaf of $E(-nN_X)$ has slope less than $-\beta_X$
and so $H^0(E(-nN_X))=0$.
Choosing a regular section of $\cO_X(nN_X)$ vanishing on a divisor $D$,
we deduce from the short exact sequence
$$
0\to E(-nN_X) \to E \to E\rest{D} \to 0
$$
that
$$
\dim H^0(E) \leq \dim H^0(E\rest{D})
= \mult(E)nN_X \leq \rk E(b+\beta_X+\mX N_X).
$$
\endproof

\para(2.3)
Now, and for the rest of the paper, let $\al$ be a positive rational
number.
Requiring $\al$ to be rational rather than real is usually not significant,
but without the restriction that $\al$ be positive many things fail,
e.g. Lemma 2.5 and its Corollaries,
and the ampleness of $\cO(p_\al,q_\al,m)$ in Section 3.2.

\say Definition 2.3.1.
A \sl Brill-Noether pair \rm on $X$ consists of a torsion-free sheaf
$E$ on $X$ together with a subspace $\LmE\sub H^0(E)$.
A \sl subpair \rm $(F,\LmF)$ of $(E,\LmE)$ consists of a
subsheaf $F\sub E$ together with a subspace $\LmF\sub\LmE\cap H^0(F)$.
A subpair is \sl proper \rm unless $F=0$ and $\LmF=0$, or $F=E$ and
$\LmF=\LmE$.
We shall call the triple $(\rk E,\dg E,\dim\LmE)$ the
\sl type \rm of the pair $(E,\LmE)$.
\endsay

\say Definition 2.3.2.
Let $(E,\LmE)$ be a Brill-Noether pair.
For any $\al\in\RR$,
the \sl $\al$-slope \rm of $(E,\LmE)$ is defined to be
$$
\mu_\al(E,\LmE)={\dg E + \al\dim\LmE \over \rk E}
$$
A Brill-Noether pair $(E,\LmE)$ is \sl $\al$-semistable \rm
if, for all proper subpairs $(F,\LmF)$,
$$
\mu_\al(F,\LmF) \leq \mu_\al(E,\LmE)
$$
Such a pair is \sl $\al$-stable \rm if
the inequality is always strict.
\endsay

By convention, the pair $(0,0)$ is $\al$-semistable
of arbitrary slope.

\para(2.4)
Notice that the restriction $\al>0$ is no great loss, because
if $\al<0$ the subpair $(E,0)\sub(E,\LmE)$ prevents any Brill-Noether
pair from being $\al$-semistable unless $\LmE=0$, in which case $\al$
is irrelevant and the definition reduces to the usual one for the
sheaf $E$.
The borderline case when $\al=0$ and $\dim\LmE>0$ has various
peculiarities which place it beyond the scope of this paper.

Even when $\al>0$, there are still restrictions on the existence of
$\al$-semistable pairs arising from the subsheaf of $E$ generated by $\LmE$.

\proclaim Lemma 2.4.
Suppose there exists an $\al$-semistable Brill-Noether pair of
type $(r,d,l)$.
\item{i)} If $l>0$, then ${d\over r}\geq -\beta_X$.
\item{ii)} If $0<l<r$, then
$\al \bigl( 1 - {l\over r} \bigr) \leq {d\over r} + \beta_X$.
Further, if there exists an $\al$-stable pair of this type,
then the inequality is strict.

\say Proof of (i).
Let $F\neq0$ be the subsheaf of $E$ generated by $\LmE$,
so that $F$ is generated by its sections and $\mu(F)\geq-\beta_X$.
Then either $F=E$ or $(F,\LmE)$ is a proper subpair and,
since $\rk F\leq \rk E$, the semistability condition implies that
$\mu(E)\geq\mu(F)$.
\endproof

\say Proof of (ii).
As for (i), let $F$ be the subsheaf generated by $\LmE$.
In this case, $(F,\LmE)$ is necessarily a proper subpair
and $\rk F\leq l$.
The semistablity condition then implies that
$\mu(F) + \al \leq {d\over r} + \al{l\over r}$,
while the stability condition requires the inequality to be strict.
\endproof

\para(2.5)
By replacing a Brill-Noether pair $(E,\LmE)$ by its `evaluation map'
$\LmE\tsr\cO_X \to E$ we may regard Brill-Noether pairs as forming
a subcategory of the abelian category
$\cC$ whose objects are arbitrary (sheaf) maps
$\zeta:\LmE\tsr\cO_X \to E$ where $\LmE$ is a finite dimensional
vector space and $E$ is any sheaf.
A morphism in $\cC$ from $\zeta:\LmE\tsr \cO_X \to E$ to
$\xi:\LmF\tsr \cO_X \to F$ is given by a linear map
$f:\LmE\to\LmF$ and a sheaf map $g:E\to F$
such that the following diagram commutes.
$$\matrix{
\LmE\tsr \cO_X & \LRA{\zeta} & E \cr
f\tsr 1 \Big\downarrow \phantom{f\tsr 1} &
& \phantom{g} \Big\downarrow g \cr
\LmF\tsr \cO_X & \LRA{\xi} & F
}$$
The `type' of an object $\zeta:\LmE\tsr\cO_X \to E$ in $\cC$
is the triple $(\rk E,\dg E,\dim\LmE)$ as for a pair
and such an object is (the evaluation map of)
a Brill-Noether pair if and only if
$E$ is torsion-free and the induced map
$H^0(\zeta):\Lam\to H^0(E)$ is injective.
Note that the subcategory so defined is not abelian
and would still not be even if we removed the restriction
that $E$ be torsion-free (see \cite{8} \S 4.1).

\say Definition 2.5.
A \sl family of objects in $\cC$ parametrised by a scheme $T$ \rm
consists of a sheaf $E_T$ on $X\times T$ flat over $T$,
a locally free sheaf $\LmE_T$ on $T$ and a map
$\zeta_T:\pi^*\LmE_T\to E_T$, where $\pi:X\times T\to T$ is
the projection map.
Another such family $(F_T,\LmF_T,\xi_T)$ is \sl isomorphic \rm
to $(E_T,\LmE_T,\zeta_T)$ if there is an invertible sheaf
$\cL$ on $T$ such that $F_T\isom E_T\tsr\pi^*\cL$,
$\LmF_T\isom\LmE_T\tsr\cL$, and under these isomorphisms
$\xi_T=\zeta_T\tsr\pi^*\cL$.
\endsay

Observe that such families pull back naturally along morphisms of
schemes and also that a family over a point is simply an
object of $\cC$.
Thus a `family of Brill-Noether pairs' is simply a family $\zeta_T$ of
objects in $\cC$ such that for each $t\in T$ the `fibre' $\zeta_t$
is a Brill-Noether pair.
A `family of $\al$-semistable Brill-Noether pairs' is defined likewise.

Now we may extend the notions of $\al$-semistability and $\al$-stability of
Brill-Noether pairs to $\cC$ by saying that an object
$\zeta:\Lam\tsr\cO_X\to E$ of type $(r,d,l)$ with $r>0$ is $\al$-semistable
if, for all subobjects $\zeta':\Lam'\tsr\cO_X\to E'$
$$
\bigl( d + \al l \bigr) \rk E'
- r \bigl( \deg E' + \al \dim\Lam' \bigr) \geq 0,
$$
and that $\zeta$ is $\al$-stable if the inequality is strict for
all proper subobjects.
By convention, no object with $r=0$ is $\al$-semistable except
0, i.e. the object with $\Lam=0$ and $E=0$.

The crucial observation is that by enlarging our category in this way,
we do not introduce any new semistable objects.

\proclaim Lemma 2.5.
An object $\zeta$ of $\cC$ of type $(r,d,l)$
is $\al$-semistable (resp. $\al$-stable)
if and only if it is (the evaluation map of) an $\al$-semistable
(resp. $\al$-stable) Brill-Noether pair.

\proof
First note that, if $\zeta:\Lam\tsr\cO_X\to E$ is
a Brill-Noether pair, then all its subobjects in $\cC$ are also.
In addition, for the non-proper subobjects of $\zeta$ the left
hand side of the above inequality is zero.
Thus all $\al$-semistable pairs are $\al$-semistable in $\cC$
and similarly all $\al$-stable pairs are $\al$-stable in $\cC$.

Conversely, any object $\zeta$ in $\cC$ has a subobject with
$\Lam'=\ker H^0(\zeta)$ and $E'=0$, and also a subobject with
$\Lam'=0$ and $E'$ the torsion subsheaf of $E$.
Both these subobjects have $\rk E'=0$ and $\deg E'\geq0$.
Hence, if $\zeta$ is $\al$-semistable and $r>0$, both
these subobjects must be 0, i.e. $\zeta$ is a Brill-Noether pair.
It is then semistable/stable as a pair if
it is semistable/stable as an object of $\cC$,
because all subpairs are certainly subobjects in $\cC$.

The conventions guarantee that the result is still true when $r=0$.
\endproof

\proclaim Corollary 2.5.1.
The full subcategory of $\cC$ consisting of $\al$-semistable
Brill-Noether pairs with fixed $\al$-slope is a noetherian
and artinian abelian category
(and hence the Jordan-H\"older Theorem holds).
The simple objects are precisely the $\al$-stable pairs.

\proof
By Lemma 2.5, this subcategory is in fact the full subcategory of
all $\al$-semistable objects in $\cC$ of type $(r,d,l)$
with $d+\al l=\mu r$, for a fixed $\mu$.
It is clear that this contains all the images, kernels
and cokernels of its maps, as well as direct sums of its objects.
Thus it is an abelian subcategory.
It is noetherian because $\cC$ is, and it is artinian because
all objects involve torsion-free sheaves and so proper subobjects must have
strictly smaller rank, thus preventing infinite descending chains.
The fact that the simple objects are the $\al$-stable pairs is immediate.
\endproof

It will also be convenient to have another variant of
the stability criterion for Brill-Noether pairs.
For this, let $P(n) = r \chi\bigl(\cO_X(n)\bigr) + d$
for any integer $n$.
Later on, we will need $n$ to be large enough that $P(n)>0$,
but this is not strictly necessary for the following.

\proclaim Corollary 2.5.2.
A Brill-Noether pair $(E,\LmE)$ of type $(r,d,l)$ with $r>0$
is $\al$-semistable if and only if, for every subsheaf $F$ of $E$
$$
\bigl( P(n) + \al l \bigr) \rk F
- r \bigl( \chi(F(n)) + \al \dim\LmE\cap H^0(F) \bigr) \geq 0
$$
Furthermore, $(E,\LmE)$ is $\al$-stable if and only if,
for every proper subsheaf, the inequality is strict.

\proof
This follows from Lemma 2.5, after observing that
$$
P(n)\rk F - r\chi(F(n)) = d\rk F - r\dg F
$$
and that every subobject of $(E,\LmE)$ in $\cC$ has the form
$(F,\LmF)$ with $\LmF\sub\LmE\cap H^0(F)$.
\endproof

\para(2.6)
A key issue in the construction of moduli spaces of the sort under
consideration is the boundedness of various families.
For this we use the following standard result
(see for example \cite{7} Lemma 8).

\proclaim Proposition 2.6.
For fixed $r,d,b$, there is a bounded family containing all torsion-free
sheaves $E$ on $X$ with $\rk E = r$, $\deg E=d$ and such that
all nonzero subsheaves $F$ of $E$ have slope  $\mu(F)\leq b$.

{}From this we deduce the following corollaries

\proclaim Corollary 2.6.1.
The family of sheaves $E$ occurring in $\al$-semistable pairs $(E,\LmE)$,
of a fixed type, is bounded.

\proof
Fixing the type fixes $\rk E$ and $\dg E$, while the semistability condition
implies that, for any nonzero subsheaf
$\mu(F)\leq {d+\al l\over r}$.
\endproof

\proclaim Corollary 2.6.2.
For fixed $r,b,k$, there is an $n_0$ such that the torsion-free sheaves
$E$, with the following properties, lie in a bounded family:
(i) $\rk E\leq r$, (ii) all nonzero subsheaves $F$ of $E$ have
slope $\mu(F)\leq b$,
and (iii) for some $n\geq n_0$,
$$
\dim H^0\bigl(E(n)\bigr)\geq \rk E\bigl( \chi\bigl(\cO_X(n)\bigr)+k \bigr).
$$

\proof
To apply the Proposition, we must obtain a lower bound on
$\mu(E)$, since, with the bounds from (i) and (ii),
there will then be only finitely many possibilities for $\rk E$ and $\deg E$.

Choose $c\leq0$, such that
${c\over r} +\max(b,0) +\beta_X +\mX N_X -\chi(\cO_X) < k$,
and $n_0$, such that $b +\mX n_0 +\beta_X \geq 0$ and
$c +\mX n_0 +\beta_X \geq 0$.
Let $E''$ be a non-zero quotient of $E$ of minimal slope, so that
$\mu(E)\geq \mu(E'')$,
and $E'$ be the kernel of the quotient map.
So $E'$ is torsion-free, but may be zero, and $E''$ is torsion-free and
semistable.
Suppose that $\mu(E'')\leq c$.
Then, for all $n\geq n_0$, we can use Lemma 2.2 to deduce
$$\eqalign{
\dim H^0(E(n))
&\leq \dim H^0\bigl(E'(n)\bigr) +\dim H^0\bigl(E''(n)\bigr) \cr
&\leq \rk E'(b +\mX n +\beta_X +\mX N_X)
     +\rk E''(c +\mX n +\beta_X +\mX N_X) \cr
&\leq \rk E (\chi(\cO_X(n)) -\chi(\cO_X) +\beta_X +\mX N_X
     +\max(b,0) +{c\over r})
}$$
which contradicts Condition (iii).
Hence we must have $\mu(E)\geq c$, as required.
\endproof

\para(2.7)
Most moduli space constructions proceed by embedding
the problem into a standard Geometric Invariant Theory (GIT) problem.
The one which arises in this paper is the natural action of $SL(V)$ on
the product of Grassmannians
$$
\grass(V\tsr H,k)\times\grass(l,V)
\eqno{(*)}
$$
where the first Grassmannian parametrises $k$-dimensional quotients of
$V\tsr H$ ($V$ and $H$ are finite dimensional vector spaces), while
the second parametrises $l$-dimensional subspaces of $V$.
This action has a (unique) linearisation on each of the line bundles
$\cO(p,q)=\pi_1^*\cO_1(p)\tsr\pi_2^*\cO_2(q)$,
where $\pi_i$ is the projection onto the $i$th Grassmannian factor
and $\cO_i(p)$ is the $p$th tensor power of the primitive ample
line bundle on that factor.
Because stability does not vary on taking tensor powers of the
linearisation, we will allow $p$ and $q$ to be rational
--- they must of course both be positive so that $\cO(p,q)$ is ample.

\proclaim Lemma 2.7.
Let $\kap_1:V\tsr H\to \Lm_1$ and $\kap_2:\Lm_2\into V$ represent a point
in the product (*).
Then $(\kap_1,\kap_2)$ is semistable with respect to the linearisation of
$SL(V)$ on $\cO(p,q)$ if and only if, for every proper subspace $U\sub V$
$$
p\th_1(U)+q\th_2(U) \geq0
$$
where
$$\eqalign{
\th_1(U) &= \dim\kap_1(U\tsr H) - {k\over \dim V}\dim U \cr
\th_2(U) &= {l\over\dim V} \dim U - \dim \Lm_2\cap U.
}$$
Further, $(\kap_1,\kap_2)$ is stable if and only if this inequality
is always strict.

\proof
This is a standard type of argument, following \cite{10} Ch.4 \S4
or \cite{11} Ch.4 \S6.
In this context, a similar argument is also given in \cite{12} Prop.5.1.1.
We give the argument as if both Grassmannians parametrised subspaces.
The reader will readily supply the modification for quotients.

Let $\kap:\Lm\into V\tsr H$ represent a point of $\grass(l,V\tsr H)$
and let $\lam:GL(1)\to SL(V)$ be a 1-parameter subgroup.
Then Mumford's numerical function, with respect to the primitive ample
line bundle, is given by
$$\eqalign{
\mu(\kap,\lam) &= \sum_{i\in\ZZ} (-i)
\bigl[ \dim\Lm\cap(V_{(\geq i)}\tsr H) -
\dim\Lm\cap(V_{(\geq i+1)}\tsr H) \bigr] \cr
&= \sum_{i\in\ZZ} \th(V_{(\geq i)})
}$$
where
$$
\th(U)= {l\over \dim V}\dim U - \dim\Lm\cap(U\tsr H).
$$
and $V_{(\geq i)}$ is the subspace of $V$ on which $\lam$ acts with weight
$\geq i$. Note that, because $\lam$ is a 1-parameter subgroup of $SL(V)$,
$$
\sum_{i\in\ZZ} i(\dim V_{(\geq i)} - \dim V_{(\geq i+1)}) = 0.
$$
Note further that for every subspace $U\sub V$
we can find a 1-parameter subgroup $\lam$
which acts with weight $(\dim V-\dim U)$ on $U$ and
weight $(-\dim U)$ on a complement of $U$.
In this case, $\mu(\kap,\lam)=\th(U) \dim V$.

The lemma then follows from the fact that $\mu$ is multiplicative under
taking tensor powers of the linearisation and additive under taking products
of spaces (provided the action and
linearisation are the obvious diagonal ones).
\endproof

%%%
%%%
%%%

\beginsection 3. Construction of the Moduli Space.

The proof of Theorem 1 follows the familiar pattern.
We construct the moduli space $G_\al(r,d,l)$ as the GIT quotient
of a projective scheme $\cR$ by a reductive group $SL(V)$.
This guarantees that $G_\al(r,d,l)$ is a projective scheme.
We show that the $SL(V)$-orbits in the open subset $\Rss\sub\cR$
of GIT semistable points are in one-one correspondence with isomorphism
classes of $\al$-semistable Brill-Noether pairs (Theorem 3.3)
amd also that GIT equivalent orbits correspond to S-equivalent pairs
(Corollary 3.4).
Finally the fact that $G_\al(r,d,l)$ is a coarse moduli space
follows in the usual way from the fact that the natural
family over $\Rss$ has the local universal property (Proposition 3.5).

Fix a possible type $(r,d,l)$ with $r>0$ and $l>0$,
and a rational number $\al>0$.
Let $P(n)=r\chi(\cO_X(n))+d$.
For all $n$ large enough that $P(n)>0$ and $\cO_X(n)$ is very ample,
let $V_n$ be a vector space of dimension $P(n)$
and $\phi_n:\cO_X\to\cO_X(n)$ be a regular section vanishing on a divisor
$D_n$.
Thus, for any torsion-free sheaf $E$,
the map $\phi_n':H^0(E)\to H^0(E(n))$, given by multiplication by $\phi_n$,
is an isomorphism of $H^0(E)$ with the kernel of the restriction map
$\rho:H^0(E(n)) \to H^0\left( E(n)\tsr\cO_{D_n} \right)$.

\para(3.1)
We write $\cQ_n$ for the quot scheme $\quot(V_n\tsr\cO_X(-n),P)$,
which parametrises quotients
$\eps:V_n\tsr\cO_X(-n)\to E$ such that the Hilbert polynomial of $E$ is $P$.
A point $\eps\in\cQ_n$ is called a `good' point if the induced map
$$
\eps'=H^0(\eps(n)):V_n\to H^0(E(n))
$$
is an isomorphism.
We write $\cG_n$ for the Grassmannian $\grass(l,V_n)$
of $l$-dimensional subspaces $\Lm\sub V_n$.
Both of these spaces carry an obvious action of $SL(V_n)$ and
two `good' points in $\cQ_n$ determine the same sheaf $E$ if and only if
they are in the same $SL(V_n)$ orbit.

Let $\cR_n$ be the closed subscheme of $\cQ_n\times\cG_n$ determined by the
condition $\Lm\sub\ker\rho\eps'$.
The action of $SL(V_n)$ clearly preserves $\cR_n$.
Furthermore, if $\eps$ is a `good' point of $\cQ_n$ and $E$ is torsion-free,
then $\bigl( E,(\phi_n')^{-1}\eps'(\Lm) \bigr)$ is a Brill-Noether
pair of type $(r,d,l)$, and two such pairs are isomorphic if and only if
the two points in $\cR_n$ are in the same $SL(V_n)$ orbit.

Now, the family of sheaves $E$ occurring in
$\al$-semistable Brill-Noether pairs $(E,\LmE)$
of type $(r,d,l)$ is bounded (Corollary 2.6.1).
Hence, we can choose $N$ large enough so that,
for all $n\geq N$, we have
$H^1(E(n))=0$ and $E(n)$ generated by its sections.
Thus any such $E$ can be represented by a point $\eps\in\cQ_n$
and, moreover, this is a `good' point.
Hence, setting $\Lm=(\eps')^{-1}\phi_n'(\LmE)$, every $\al$-semistable pair
$(E,\LmE)$ corresponds to a unique $SL(V_n)$ orbit in $\cR_n$.

\para(3.2)
We now describe the GIT problem on $\cR_n$, for which
the stability conditions on orbits will be seen to correspond exactly
to the stability conditions on pairs.
We follow Simpson's method, which is based on Grothendieck's
original construction of the quot scheme.

We may choose $M_n$ sufficiently large that, for $m\geq M_n$,
the functor, which takes $\eps:V_n\tsr\cO_X(-n)\to E$ to
$$
\eps''=H^0(\eps(n+m)):V_n\tsr H^0(\cO_X(m))\to H^0(E(n+m)),
$$
determines a closed embedding of
$\cQ_n$ in $\grass(V_n\tsr H^0(\cO_X(m)),P(n+m))$.
Thus we obtain a closed embedding
$$
\cQ_n\times\cG_n\into\grass(V_n\tsr H^0(\cO_X(m)),P(n+m))\times\grass(l,V)
$$
Let $\cO(p,q,m)$ be the pullback to $\cQ_n\times\cG_n$ of the line bundle
$\cO(p,\mX mq)$ on the product of Grassmannians (c.f. 2.7).

For any subspace $U\sub V_n$, let $F_U=\eps(U\tsr\cO_X(-n))$ and
$K_U=\ker\eps\rest{U\tsr\cO_X(-n)}$.
Because $\eps$ and $U$ lie in a bounded family, $F_U$ and $K_U$
do also.
Hence, we may further suppose that $M_n$ is sufficiently large that,
for all $m\geq M_n$,
$$
H^1(F_U(n+m))=H^1(K_U(n+m))=0.
$$
Then $\eps''(U\tsr H^0(\cO_X(m))=H^0(F_U(n+m))$, which is a vector
space of dimension $\chi(F_U(n+m))$.
Applying Lemma 2.7, with $V=V_n$, $H=H^0(\cO_X(m))$ and $k=P(n+m)$,
yields a stability condition for a point in $\cQ_n\times\cG_n$ in
terms of
$$\eqalign{
\th_1(U) &= \chi(F_U(n+m))-{P(n+m)\over P(n)}\dim U \cr
&= \mX{m\over P(n)} \bigl( P(n)\rk F_U - r\dim U \bigr)
+ \chi(F_U(n))-\dim U
}$$
and
$$
\th_2(U) = {l\over P(n)} \dim U - \dim \Lm\cap U.
$$
Indeed, if $\al$ is a positive rational number and we set
$p_\al = P(n)+\al l$ and $q_\al=\al r$,
then stability with respect to the linearisation $\cO(p_\al,q_\al,m)$
depends on
$$
p_\al\th_1(U) + \mX mq_\al\th_2(U)
= \mX m\th_\al(U)
+ \bigl( P(n)+\al l \bigr) \bigl( \chi(F_U(n)) - \dim U \bigr)
$$
where
$$
\th_\al(U) =
(P(n)+\al l) \rk F_U - r \bigl( \dim U + \al \dim\Lm\cap U \bigr).
$$
This quantity should be compared with the one which occurs in Corollary 2.5.2.

Now, we will say that a point in $\cQ_n\times\cG_n$ is
$(p,q)$-semistable (resp. stable) if it is GIT semistable (resp. stable)
with respect to the linearisation $\cO(p,q,m)$ for all $m\gg0$.
One can readily see that this condition can be formulated without $m$.

\proclaim Lemma 3.2.
A point $(\eps,\Lm)\in\cQ_n\times\cG_n$ is $(p_\al,q_\al)$-semistable
if and only if, for all subspaces $U\sub V_n$,
we have $\th_\al(U) \geq 0$ and, in the case of equality,
$\chi(F_U(n)) \geq \dim U$.
Further, $(\eps,\Lm)$ is $(p_\al,q_\al)$-stable
if and only if, for all proper subspaces,
one of the inequalities is strict.

\proof
We simply need to observe that the quantities $\dim U$ and $\chi(F_U(n))$
are bounded above and below, independent of $m$, because $U$ and $F_U$ lie
in bounded families.
Hence we can choose $m$ large enough that
${1\over m}|\chi(F_U(n))-\dim U|(P(n)+\al l)$
is less than the smallest possible
positive value of $\mX\th_\al(U)$,
which depends only on the denominators of $\al$ and $P(n)$,
since $\mX\rk F_U$ and $\mX r$ are both integers (see 2.1).
\endproof

\para(3.3)
We are now in a position to prove that notions of
$\al$-stability of Brill-Noether pairs and of $(p_\al,q_\al)$-stability
of points in $\cR_n$ correspond exactly.

\proclaim Theorem 3.3.
There is an integer $N$ such that, for all $n\geq N$, the following hold
\item{i)} If $(E,\LmE)$ is an $\al$-semistable (resp. $\al$-stable)
Brill-Noether pair, then the corresponding orbit in
$\cR_n$ is $(p_\al,q_\al)$-semistable (resp. $(p_\al,q_\al)$-stable).
\item{ii)} If $(\eps,\Lm)\in\cQ_n\times\cG_n$ is
$(p_\al,q_\al)$-semistable, then $H^1(E(n))=0$, and
$\eps$ is a `good' point of $\cQ_n$, and $E$ is torsion-free.
\item{iii)} If $(\eps,\Lm)\in\cR_n$ is $(p_\al,q_\al)$-semistable
(resp. $(p_\al,q_\al)$-stable), then the associated Brill-Noether
pair is $\al$-semistable (resp. $\al$-stable).

\say Proof of (i).
Let $(E,\LmE)$ be $\al$-semistable and $(\eps,\Lm)$ a point in the
corresponding orbit in $\cQ_n\times\cG_n$.
Let $U$ be any proper subspace of $V_n$
and $F_U=\eps(U\tsr\cO_X(-n))$, as in the previous section.
Let $F$ be the smallest subsheaf of $E$ which contains $F_U$ and for which
$E/F$ is torsion-free.
Note that $\rk F=\rk F_U$.
Since the induced map $\eps':V_n\to H^0(E(n))$ is an isomorphism,
$U$ can be regarded as a subspace of both $H^0(F_U(n))$ and $H^0(F(n))$.
Further, all the sections in $\Lm\sub H^0(E(n))$
vanish on $D_n$ and so $\Lm\cap U$ can be identified with a subspace of
$\LmE\cap H^0(F)$. (Note that this requires $E/F$ to be torsion-free.)

Now, by Corollary 2.6.2, we can choose $N$ large enough so that,
for all $n\geq N$, either $H^1(F_U(n))=0$, or
$$
\dim H^0\bigl( F_U(n) \bigr) + \al l
< \rk F_U \left( {P(n)+\al l \over r} \right).
$$
In the second case, because $\dim U \leq \dim H^0(F_U(n))$ and
$\dim\Lm\cap U \leq l$, we deduce that $\th_\al(U)>0$.
In the first case, we have
$$
\chi(F_U(n)) =\dim H^0(F_U(n)) \geq\dim U
$$
and also $H^1(F(n))=0$, so that $\chi(F(n)) \geq\dim U$ and
$$
\th_\al(U) \geq
\bigl( P(n) + \al l \bigr) \rk F
- r \bigl( \chi(F(n)) + \al \dim\LmE\cap H^0(F) \bigr)
$$
Thus, by Corollary 2.5.2 and Lemma 3.2, if $(E,\LmE)$ is $\al$-semistable,
then $(\eps,\Lm)$ is $(p_\al,q_\al)$-semistable.
Similarly if $(E,\LmE)$ is $\al$-stable and $F\neq E$,
then $(\eps,\Lm)$ is $(p_\al,q_\al)$-stable.
On the other hand, if $F=E$, then $\chi(F(n))=\dim V_n >\dim U$
and the above inequality on $\th_\al(U)$ is strict with right-hand-side zero.
Hence, $(\eps,\Lm)$ is $(p_\al,q_\al)$-stable in this case as well.
\endproof

\say Proof of (ii).
By Serre duality, $H^1(E(n))=0$ if and only if
$\Hom(E(n),\om_X)=0$, where $\om_X$ is the dualising sheaf
(\cite{5} III.7).
So, suppose we have a map $\psi:E(n)\to\om_X$ and let
$U$ be the kernel of the map
$$
H^0(\psi)\eps':V_n\to H^0(\om_X).
$$
Observe that $\dim U\geq P(n)-g$,
where $g=\dim H^0(\om_X)$,
and so the semistability condition implies that
$$\eqalign{
0 &\leq
(P(n)+\al l)\rk F_U - r(P(n)-g) \cr
&= r(g + \al l) - (P(n)+\al l)(r-\rk F_U)\cr
}$$
By choosing $N$ sufficiently large that, for all $n\geq N$,
we have $P(n)+\al l > \mX r(g+\al l)$, we guarantee that $\rk F_U=r$.
But $F_U(n)\sub\ker\psi$, and so $\rk(\ker\psi)=r$ and $\rk(\im\psi)=0$.
But $\om_X$ is torsion-free and so $\psi=0$, as required.

Thus $V_n$ and $H^0(E(n))$ have the same dimension and, to show that
$\eps$ is a good point of $\cQ_n$, it is sufficient to show that $\eps'$
is injective.
So, let $U=\ker\eps'$.
Then $F_U=0$ and the semistability condition reduces to
$$
0\leq -r(\dim U + \al\dim\Lm\cap U)
$$
which implies that $U=0$, as required.

Finally, to show that $E$ is torsion-free, let $T\sub E$ be a torsion
subsheaf and observe that, since $V_n\isom H^0(E(n))$, we can choose
$U\isom H^0(T(n))$.
But then $F_U\sub T$ and so $\rk F_U=0$.
Just as above, this implies that $U=0$ and hence that $T=0$.
\endproof

\say Proof of (iii).
We use the criterion in Corollary 2.5.2, for the semistability/stability
of Brill-Noether pairs.
Let $(\eps,\Lm)$ be a $(p_\al,q_\al)$-semistable point in $\cR_n$
and $(E,\LmE)$ the associated Brill-Noether pair.
Let $F\sub E$ be any subsheaf.
Because $V_n\isom H^0(E(n))$,
$F$ determines a subspace $U\isom H^0(F(n))$.
Because $F$ is torsion-free, $\LmE\cap H^0(F)$ can be identified
with a subspace of $\Lm\cap U$.
Hence $\dim \Lm\cap U \geq\dim \LmE\cap H^0(F)$, while also
$\dim U\geq\chi(F(n))$, and $\rk F_U\leq\rk F$, because $F_U\sub F$.
Therefore,
$$
\bigl( P(n) + \al l \bigr) \rk F
- r \bigl( \chi(F(n)) +\al \dim \LmE\cap H^0(F) \bigr)
\geq \th_\al(U) \geq0
$$
and hence, $(E,\LmE)$ is $\al$-semistable.

We must now show that, if $(E,\LmE)$ is strictly $\al$-semistable
(i.e. is not $\al$-stable), then $(\eps,\Lm)$ is
strictly $(p_\al,q_\al)$-semistable.
So, suppose that $F$ is a proper subsheaf with
$$
\bigl( P(n) + \al l \bigr) \rk F
- r \bigl( \chi(F(n)) +\al \dim \LmE\cap H^0(F) \bigr) =0
$$
and hence $\th_\al(U)=0$, with $U\isom H^0(F(n))$, as above.
Now, $(F,\LmE\cap H^0(F))$ must be an $\al$-semistable pair of
one of finitely many types, and so, by Corollary 2.6.1, $F$ lies in
a bounded family.
We may therefore choose $N$ sufficiently large that, for $n\geq N$,
$H^1(F(n))=0$ and $F(n)$ is generated by its sections.
Hence, $U$ must be a proper subspace and $F_U=F$,
so $\chi(F_U(n))=\dim U$.
Thus, by Lemma 3.2, $(\eps,\Lm)$ is not $(p_\al,q_\al)$-stable.
\endproof

\para(3.4)
{}From now on, we drop the subscript $n$ on $V_n$, $\cR_n$, etc. and
implicitly assume that $n$ has been chosen large enough that all the
properties required in Theorem 3.3 hold.

We now show that GIT equivalence in $\Rss$ (the set of semistable
points in $\cR$) corresponds to S-equivalence of Brill-Noether pairs.
This happens because there is a sufficiently good correspondence between
1-parameter subgroups of $SL(V)$ and filtrations of pairs.

\proclaim Proposition 3.4.
Let $(\eps,\Lm)$ be a point in $\Rss$ corresponding to an
$\al$-semistable pair $(E,\LmE)$.
\item{i)} If $\lam:GL(1)\to SL(V)$ is a 1-parameter subgroup for which
the limit
$$
(\eps_0,\Lm_0)=\lim_{t\to0} \lam(t)(\eps,\Lm)
$$
is also in $\Rss$, then $\lam$ induces a filtration of $(E,\LmE)$
by $\al$-semistable subpairs of the same $\al$-slope
and $(\eps_0,\Lm_0)$ corresponds to
the associated graded pair for this filtration.
\item{ii)} Given any filtration of $(E,\LmE)$ by $\al$-semistable subpairs
of the same $\al$-slope,
there is a 1-parameter subgroup $\lam$ which induces this filtration
as in (i).

\say Proof of (i).
First note (following \cite{15} p.33)
that $\lam$ induces a $\ZZ$-parametrised filtration
$V_{(\geq i)}$ of $V$, consisting of the subspaces on which $\lam$ acts
with weight at least $i$, together with an isomorphism $V\isom \gr V$,
where
$$
\gr V = \bigoplus_{i\in\ZZ} V_{(\geq i)}/V_{(\geq i+1)}.
$$
This filtration in turn generates filtrations
$$\eqalign{
E_{(\geq i)} &= \eps\bigl( V_{(\geq i)}\tsr\cO_X(-n) \bigr) \cr
\Lm_{(\geq i)} &= \Lm\cap V_{(\geq i)}
}$$
of $E$ and $\Lm$ respectively, and, consequently, associated graded
objects $\gr\eps:V\tsr\cO_X(-n) \to \gr E$ and $\gr\Lm\sub V$.

Now, $\eps_0=\gr\eps$ and $\Lm_0=\gr\Lm$ (\cite{15} Lemma 1.26) and
hence, since $(\eps_0,\Lm_0)\in\Rss$, Theorem 3.3(ii) shows,
firstly, that $\gr E$ is torsion-free and hence, inductively,
that $E/E_{(\geq i)}$ is torsion-free, and secondly, that $\gr\eps$ is
a `good' point and so $V_{(\geq i)}\isom H^0(E_{(\geq i)}(n))$.
Thus (c.f. proof of Theorem 3.3(i) ), $\Lm_{(\geq i)}$ can be identified
with $\LmE\cap H^0(E_{(\geq i)})$.
This subspace, together with $E_{(\geq i)}$ itself, make up a
subpair $(E,\LmE)_{(\geq i)}$ in a $\ZZ$-filtration of $(E,\LmE)$,
whose associated graded pair corresponds to $(\eps_0,\Lm_0)$.
This associated graded pair is $\al$-semistable with $\al$-slope $\mu$, say.
Hence each summand has $\al$-slope $\mu$ and, inductively, each subpair
$(E,\LmE)_{(\geq i)}$ has $\al$-slope $\mu$ and so, being a subpair
of an $\al$-semistable pair of $\al$-slope $\mu$, is itself $\al$-semistable.
\endproof

\say Proof of (ii).
Given any $\ZZ$-filtration of $(E,\LmE)$ by $\al$-semistable subpairs of the
same $\al$-slope, the subspace of sections associated to $E_{(\geq i)}$
must be $\LmE\cap H^0(E_{(\geq i)})$.
Furthermore, $E_{(\geq i)}(n)$ is generated by its sections
(c.f. proof of Theorem 3.3(iii)).
Thus, the filtration is generated by $V_{(\geq i)}=H^0(E_{(\geq i)}(n))$,
which can in turn be induced by a 1-parameter subgroup of $GL(V)$.
By altering the indexing of the filtration, but not its terms, we may choose
this 1-parameter subgroup to be in $SL(V)$.
\endproof

\proclaim Corollary 3.4.
Two points $(\eps_1,\Lm_1)$ and $(\eps_2,\Lm_2)$ in $\Rss$ are
GIT equivalent if and only if the corresponding $\al$-semistable pairs
$(E_1,\LmE_1)$ and $(E_2,\LmE_2)$ are S-equivalent.

\proof
If two points are GIT equivalent they have a common closed orbit in their
orbit closures and hence there are one parameter subgroups whose limits
are this same closed orbit. Hence, by (i) above, the corresponding pairs
are S-equivalent.
On the other hand, by (ii) above, if the corresponding pairs are
S-equivalent, then there are 1-parameter subgroups which take the
points (in the limit) to the same orbit and so the points are GIT
equivalent.
\endproof

\para(3.5)
We conclude by showing that there is a family of $\al$-semistable
Brill-Noether pairs on $\Rss$ with the local universal property.
Below, we shall use $\pi$ to denote the projection $T\times X\to T$
for any $T$, which will be apparent from the context.
For a sheaf $F$ on $T\times X$ we shall write $F(n)$ for
$F\tsr(p_X)^*\cO_X(n)$ where $p_X:T\times X\to X$ is the other projection.

For $E_T$ on $T\times X$ and $\LmE_T$ on $T$,
we shall exploit the natural isomorphism
$$
\Hom(\pi^*\LmE_T,E_T)\isom
\Hom(\LmE_T,\pi_*E_T)
$$
to regard a family of Brill-Noether pairs as specified
by a map $j_T:\LmE_T\to\pi_*E_T$ such that,
for each $t\in T$, $E_t$ is torsion-free and the induced map
$\LmE_t\to H^0(E_t)$ is injective.
This last condition implies that $j_T$ is a sheaf injection,
and we thus recover Definition 1.11 of \cite{13}.

Now the family on $\Rss$ comes from
(i) the universal quotient sheaf $E_\cQ$ on $\cQ\times X$,
for which $\pi_*E_\cQ$ can be identified with
a subsheaf of $V\tsr\cO_\cQ$ over the good points corresponding
to torsion-free sheaves,
and (ii) the universal subsheaf $\Lm_\cG\sub V\tsr\cO_\cG$ on $\cG$.
More accurately, let $E_\Rss$ be the pull-back of $E_\cQ$ to
$\Rss\times X$ and let $\LmE_\Rss$ be the pull-back of $\Lm_\cG$
to $\Rss$.
Observe that, by definition of $\Rss$,
$\LmE_\Rss$ is a subsheaf of $\pi_*E_\Rss$
and further $\LmE_x$ is contained in $H^0(E_x)$ for all $x\in\Rss$.
Thus $\bigl( E_\Rss,\LmE_\Rss \bigr)$ is a family of Brill-Noether
pairs, which by construction are all $\al$-semistable.

\proclaim Proposition 3.5.
Let $(E_T,\LmE_T)$ be a family of $\al$-semistable Brill-Noether pairs
which is parametrised by a scheme $T$.
Then each $t\in T$ has an open neighbourhood $U\sub T$ over which
there is defined a morphism $U\to\Rss$
such that $(E_U,\LmE_U)$ is isomorphic
to the pull-back of $\bigl( E_\Rss,\LmE_\Rss \bigr)$.

\proof
As in Section 3.4, suppose that $n$ has been chosen large enough
that all the properties required in Theorem 3.3 hold.
Then, since all the pairs in the family are semistable,
$\pi_*E_T(n)$ is locally-free over $T$ of rank equal to the dimension
of $V$.
Hence, for $t\in T$, we may choose an open neighbourhood $U\sub T$
and an isomorphism of $\pi_*E_U(n)$ with $V\tsr\cO_U$.
This determines a surjective map $V\tsr\cO_{U\times X}(-n)\to E_U$ and hence
a morphism $U\to\cQ$.
We also have a sheaf inclusion
$\LmE_U\into \pi_*E_U\into V\tsr\cO_U$ which is injective on
every fibre, hence determines a morphism $U\to\cG$.
These two morphisms determine a morphism
$U\to\Rss\sub\cQ\times\cG$.
It is clear from the construction that the original
family is the pullback of the family on $\Rss$.
\endproof

Thus the categorical quotient of $\Rss$ is a coarse moduli space for
families of $\al$-semistable Brill-Noether pairs.

%%%
%%% bibliography
%%%

\beginsection References.

\item{1.} E. Arbarello, M. Cornalba, P.A. Griffiths \& J. Harris,
Geometry of Algebraic Curves, Vol. I, Springer-Verlag, New York, (1985).

\item{2.} A. Bertram,
Stable pairs and stable parabolic pairs,
Harvard preprint, 1992.

\item{3.} S.B. Bradlow \& G.D. Daskalopoulos,
Moduli of stable pairs for holomorphic bundles over Riemann surfaces,
Int. J. Math. 2 (1991) 477--513.

\item{4.} S. Bradlow, G.D. Daskolopoulos, O. Garcia-Prada
\& R. Wentworth,
Stable augmented bundles over Riemann surfaces,
to appear in {\it Vector Bundles in Algebraic Geometry},
Durham 1993.

\item{5.} R. Hartshorne,
{\it Algebraic Geometry},
Springer-Verlag, New York (1977).

\item{6.} D. Huybrechts \& M. Lehn,
Stable pairs on curves and surfaces,
J. Algebraic Geometry, to appear.

\item{7.} J. Le Potier,
Faisceaux semi-stables sur les courbes,
Jussieu notes (1991).

\item{8.} J. Le Potier,
Syst\`emes coh\'erents et structures de niveau,
Ast\'erisque 214 (1993).

\item{9.} J. Le Potier,
Faisceaux semistables et syst\`emes coh\'erents,
to appear in {\it Vector Bundles in Algebraic Geometry},
Durham 1993.

\item{10.} D. Mumford \& J. Fogarty,
Geometric Invariant Theory, 2nd Edition, Springer-Verlag, 1982.

\item{11.} P.E. Newstead,
Introduction to Moduli Problems and Orbit Spaces,
Tata Lecture Notes, Springer-Verlag, 1978.

\item{12.} M.S. Narasimhan \& G. Trautmann,
Compactification of $M_{\PP_3}(0,2)$ and Poncelet pairs of conics,
Pacific J. of Math. 145 (1990) 255--365.

\item{13.} N. Raghavendra \& P.A. Vishwanath,
Moduli of pairs and generalized theta divisors,
Tohoku Math. J., to appear.

\item{14.} C.S. Seshadri,
Fibr\'es vectoriels sur les courbes alg\'ebriques,
Ast\'erisque 96 (1982).

\item{15.} C. Simpson,
Moduli of representations of the fundamental group of a smooth
projective variety - I \& II, Toulouse preprint 17 (1992).

\item{16.} M. Thaddeus,
Stable pairs, linear systems and the Verlinde formula,
Invent. Math. 117 (1994) 317--353.

\bye